\documentclass[a4paper,11pt]{article}

\usepackage{url}
\usepackage{amssymb,amsmath,amsthm,amsfonts}
\usepackage[all]{xy}

\newenvironment{demo}{\noindent {\sl Proof}. \ }{\qed}

\newtheorem{stheoreme}{Theorem}[subsection]
\newtheorem{sdefin}[stheoreme]{Definition}
\newtheorem{sprop}[stheoreme]{Proposition}
\newtheorem{slemme}[stheoreme]{Lemma}
\newtheorem{scorol}[stheoreme]{Corollary}
\newtheorem{sremark}[stheoreme]{Remark}

\newcommand\zum[2]{\sum_{\substack{#1\\#2}}}
\newcommand\ozum[2]{\bigoplus_{\substack{#1\\#2}}}

\def\ssum{\sum\limits}

\font\tenCal=cmsy10
\newfam\Calfam
\textfont\Calfam=\tenCal

\def\disy{\displaystyle}

\def\Hom{\mathop{\hbox{\textrm Hom}}\nolimits}

\def\qed{\hfill{$\sqcap\!\!\!\!\sqcup$}}

\def\wt{\widetilde}
\def\what{\widehat}

\def\surl{\overline}

\def\L{\Lambda}

\def\a{\alpha}
\def\b{\beta}
\def\g{\gamma}

\def\CC{{\mathbb C}}
\def\blanc{\ \ \ \ \ \ \ \ \ \ }
\def\esp{\ \ \ \ \ }

\title{Twisted Poincar\'e duality for some quadratic Poisson algebras}

\author{St\'ephane Launois\thanks{This research was supported by a Marie Curie Intra-European
 Fellowship within the $6^{\mbox{th}}$ European Community Framework
 Programme} \ and Lionel Richard\thanks{Supported by  EPSRC Grant
EP/D034167/1}}

\date{}

\begin{document}

\maketitle
\begin{abstract}
We exhibit a Poisson module restoring a twisted Poincar\'e duality between Poisson homology and cohomology for the polynomial algebra $R={\mathbb C}[X_1,\ldots,X_n]$ endowed with Poisson bracket arising from a uniparametrised quantum affine space. This Poisson module is obtained as the semiclassical limit of the dualising bimodule for Hochschild homology of the corresponding quantum affine space.
As a corollary we compute the Poisson cohomology of $R$, and so retrieve a result obtained by direct methods (so completely different from ours) by Monnier.
\end{abstract}


\section{Introduction}

\hspace{1cm} Given a Poisson algebra, its Poisson cohomology
provides important information about the Poisson structure (the
Casimir elements are reflected by the degree zero cohomology,
Poisson derivations modulo Hamiltonian derivations by the degree
one,...). Computing this cohomology is in general difficult. It
has been achieved in some particular cases, see for instance
\cite{Pichereau} and references therein. One way to study a
Poisson algebra is to consider a deformation. If the deformation
is ``nice'', its properties should reflect the corresponding
properties of the original Poisson algebra. Results in this spirit
have been obtained for instance in \cite{FT}, or in the framework
of symplectic varieties in \cite{AL98} and recently in \cite{AF},
where roughly speaking the Poisson homology is shown to match the
Hochschild homology of the deformation at least in small degrees.

The aim of this paper is to illustrate this idea that quantisation
can provide some intuition on the study of Poisson (co)homology.
More precisely, we consider the polynomial algebra $R={\mathbb
C}[X_1,\ldots,X_n]$ endowed with the bracket
$\{X_i,X_j\}=a_{ij}X_iX_j$, where $(a_{ij})\in M_n({\mathbb Z})$
is a skew-symmetric matrix. The Poisson cohomology of this algebra
has been computed in \cite{Monnier} by direct computation. Another
approach in order to compute this Poisson cohomology consists in
establishing a duality between Poisson homology and cohomology. In
general such a duality does not occur. For instance for $n=2$ and
$a_{12}=1$ the cohomology spaces in degree 0,1 and 2 have
respectively dimension 1,2,2, whereas the homology space of degree
2 is null, and infinite-dimensional in degree 0 and 1. Note that
this lack of nontrivial Poisson homology in the higher degrees explains the  ``dimension drop'' appearing in the Hochschild homology of the quantum plane (see \cite{manawa}, \cite{Wa}, \cite{Si}, \cite{HK1}, \cite{BZ}, \cite{connesdbviolette} for details on dimension drop in this and other situations).

  The starting
point of this work was the following question: knowing that there
is no Poincar\'e duality between Poisson homology and cohomology
in this case (see the paper \cite{RV} for $n=2$), can it be
replaced by a ``twisted duality''? Such a twisted duality for
Poisson (co)homology appears in the literature (see for instance
\cite{Hueb99}, \cite{EvensLuWeinstein}, \cite{Xu}), in a somewhat
abstract way. In the present paper we make this twisted duality
explicit on the complexes that are used to compute the Poisson
(co)homology of the Poisson algebra under consideration.

Our method here relies on the fact that the Poisson algebra $R$ is
the semiclassical limit of a uniparametrised  quantum affine space
$U$ admitting a dualising bimodule for Hochschild (co)homology.
This bimodule is easily defined as the algebra $U$ itself, with
product twisted on the left by an automorphism. This allows us to
define a Poisson module $M$ over $R$ as the semiclassical limit of
the dualising bimodule of $U$. Then we show that the Poisson
cohomology of $R$ is dual to the Poisson homology with values in
$M$; in other words, $M$ is a ``Poisson dualising module'' in the
sense that it restores a duality between Poisson homology and
cohomology. Finally we use this twisted Poisson duality to compute
the Poisson cohomology of $R$,  retrieving so a
 result of Monnier.

 Although we work here with the semiclassical limit of quantum affine space, it is likely that this method will apply to other algebras admitting a ``twisted'' Poincar\'e duality for the Hochschild homology at the quantised level, at least when the automorphism used for this twisting can be expressed simply as multiplying generators by a power of the parameter $q$. This is the case for instance for quantised coordinate rings of semisimple complex algebraic groups, as proved in \cite{BZ}. Starting from the Van den Bergh duality, \cite{VdB}, at the quantum level, it would be very interesting to retrieve this twisted Poisson duality for semiclassical limits  thanks to spectral sequences \`a la Brylinski, \cite{brylinski}, see also \cite{maszczyk}. We plan to go back to these questions in a subsequent paper.

The plan of this paper  is as follows. First, we show that the automorphism provided by Sitarz to solve the ``dimension drop'' for the Hochschild homology of the quantum space provides a twisted duality. Next,  we produce a Poisson module providing duality between Poisson homology and cohomology of the semiclassical limit $R$. Finally, we compute the Poisson homology with values in this module.

Throughout this paper we will use the usual following notation for a monomial: if $x_1,\ldots,x_n$ are variables and $\alpha=(\alpha_1,\ldots,\alpha_n)\in{\mathbb N}^n$ then $x^{\alpha}=x_1^{\alpha_1}\ldots x_n^{\alpha_n}$. We will also denote by $\epsilon_i$ the $i$th vector of the canonical basis of ${\mathbb Z}^n$.

\section{Twisted Poincar\'e duality in Hochschild (co)homology.}

\hspace{1cm}
In the paper \cite{Si},  Sitarz provides an automorphism of  quantum affine space restoring what one should expect as the ``good'' Hochschild dimension of these algebras. Namely,  Wambst proved in \cite{Wa} that in the generic case all Hochschild homology groups are null in degree greater than 1.
 Sitarz proved in his paper that the ``twisted'' Hochschild homology, for a particular choice of twisting automorphism, provides a vector space of dimension 1 for the homology group of degree $n$ of the quantum affine $n$-space.

Re-interpreting certain results of \cite{Rhh} in terms of ``twisted'' Poincar\'e duality (see \cite{VdB}, \cite{BZ}), we show here that the automorphism determined by  Sitarz actually provides such a duality. This could be as well seen as a corollary of Van den Bergh's results, \cite[Proposition 2]{VdB}, but we make it explicit here at the level of the complexes themselves.
The method used, in terms of ``up to a sign'' commuting diagrams, is the same as the one we are going to use for Poisson homology.
The following computations are done for quantum affine space, but one can easily check that they remain valid for the ``mixed crossed algebras'' defined in \cite{Rhh}.

\subsection{Resolution} \label{ck}

\hspace{1cm}
We say that a matrix  $Q=(q_{ij})\in M_n({\mathbb C})$ is multiplicatively skew-symmetric if $q_{ij}q_{ji}=q_{ii}=1$ for all $i,j$. Let us recall here the following two quantum algebras.

\begin{sdefin} \label{defqalg}
Let $V$ be a ${\mathbb C}$-vector space of dimension $n$ with basis $(v_1,\ldots,v_n)$, and let $Q=(q_{ij})_{1\leq i,j\leq n}\in M_n({\mathbb C}^*)$ be a multiplicatively skew-symmetric matrix.
\begin{enumerate}
\item
The \emph{quantum affine space} $S_QV={\mathcal O}_Q(\CC^n)$ is the $\CC$-algebra generated by $n$ generators $v_1,\ldots,v_n$ with relations $v_iv_j-q_{ij}v_jv_i$ for all $1\leq i,j\leq n$.

 It is well-known that $S_QV$ is a noetherian domain, and admits the monomials $\{v^{\a}\}_{\a\in{\mathbb N}^n}$ as a PBW-basis.
\item
The \emph{quantum exterior algebra} $\L_QV$ is the $\CC$-algebra  generated by $n$ generators $v_1,\ldots,v_n$ with relations $v_iv_j=-q_{ij}v_jv_i$ for all $1\leq i,j\leq n$.
It admits the monomials $\{v^{\beta}\}_{\beta\in\{0,1\}^n}$ as a $\CC$-linear basis.
\end{enumerate}
\end{sdefin}
\begin{sremark}{\rm
This definition of the quantum exterior algebra differs from the one given in \cite{bg} for instance. However, it is the one provided in \cite{Wa}, giving rise to the quantum Koszul complex described below.}
\end{sremark}
Set $U=S_QV$.
It follows from \cite[Proposition 10.3]{Wa} that the following is a free resolution of $U$ as a $U^e=U\otimes U^{op}$-module.

The vector space $\L_QV$
is ${\mathbb N}$-graded, by defining
the degree of a monomial
 $v_1^{\b_1}\wedge\ldots\wedge v_n^{\b_n}$ to be  $|\b|$.
For all
 $*\in{\mathbb N}$, one denotes by $\L^*_QV$ the homogeneous subspace of $\L_QV$ of degree
  $*$. Then  the space
 $U\otimes \L_QV\otimes U$ is graded by  the degree of
 $\L_QV$, and $U\otimes \L_QV\otimes U$ becomes a differential complex, with the differential
 $\partial$ defined for all
 $a,b\in U$ by:
\begin{equation} \label{ddd}
\begin{array}{l}
\partial(a\otimes v_{i_1}\wedge\ldots\wedge v_{i_*} \otimes b)=\\
\blanc {\disy \sum_{k=1}^* (-1)^{k-1}\bigg(\Big( \prod_{s<k}q_{i_s,i_k}\Big)av_{i_k}\otimes v_{i_1}\wedge\ldots\wedge \what v_{i_k} \wedge\ldots \wedge v_{i_*}\otimes b}  \\
 \blanc \blanc  -{\disy  \Big(\prod_{s>k}q_{i_k,i_s}\Big)a\otimes v_{i_1}\wedge\ldots \wedge\what v_{i_k}\wedge \ldots \wedge v_{i_*}\otimes v_{i_k} b} \bigg)
\end{array}
\end{equation}
where $v_{i_1}\wedge\ldots \wedge \what v_{i_k}\wedge \ldots \wedge v_{i_*}$ is the outer product $v_{i_1}\wedge\ldots \wedge v_{i_*}$ without $v_{i_k}$.

\subsection{Hochschild homology, cohomology}

\hspace{1cm}
We recall some of the homological framework that is used in \cite{Wa}, \cite{Si}.

\begin{sdefin}
Let ${\mathcal A}$ be a $\CC$-algebra. An automorphism $\sigma$ of ${\mathcal A}$ is said to be a \emph{scaling automorphism} if there exists a basis $\{a_i\}_{i\in I}$ of ${\mathcal A}$ as a vector space such that $\sigma(a_i)=p_ia_i$ with $p_i\in \CC^*$ for all $i\in I$.
\end{sdefin}

\begin{sdefin} \label{defhhtordue}
Consider a $\CC$-algebra $ {\mathcal A}$, $\sigma$ an automorphism of $A$, and let ${}_{\sigma}{\mathcal A}$ be the ${\mathcal A}^e$-module that is ${\mathcal A}$ as a vector space with module structure: $b.(a_0,a_1)=\sigma(a_0)ba_1$. Denote by $K({\mathcal A})$ a projective resolution of ${\mathcal A}$ by ${\mathcal A}^e$-modules.
\begin{enumerate}
\item The \emph{invariant twisted Hochschild homology} is the homology of the subcomplex of\\
 ${}_{\sigma}{\mathcal A}\otimes_{{\mathcal A}^e}K({\mathcal A})$ consisting only  of $\sigma$-invariant elements.
\item The \emph{twisted Hochschild homology} is the homology of the quotient of the (usual) Hochschild complex by the image of the map $1-\sigma$.
\end{enumerate}
\end{sdefin}

\begin{slemme}[\cite{Si}, Lemma 2.2] For any scaling automorphism $\sigma$ the corresponding invariant twisted Hochschild homology is isomorphic to the twisted Hochschild homology.
\end{slemme}

So the twisted Hochschild homology is computed as the $\sigma$-invariant part of the Hochschild homology with values in ${}_{\sigma}{\mathcal A}$.
Thus we are interested here in $H(U,{}_{\sigma}U)$.
 The following scaling automorphism for the algebra $S_QV$ is defined in \cite{Si}.

\begin{sdefin}  The automorphism $\sigma$ of $S_QV$ defined by $\sigma(v_i)=p_iv_i$, with $p_i=\prod_{j=1}^nq_{j,i}$ is a scaling automorphism, called the  \emph{canonical} scaling automorphism.
\end{sdefin}

Let us recall - with adapted notation - the results obtained in \cite[Section 6]{Wa}  and \cite[Section 3]{Si}.
The twisted
 Hochschild homology of $S_{Q}$ is given by the complex below:
$$K^{\sigma}(S_QV)=\bigoplus_{\tiny{\begin{array}{c} \a\in {\mathbb N}^n\\ \g\in\{0,1\}^n \end{array}}} \CC.v^{\a}\otimes v^{\g},$$
\begin{equation} \label{dwq}
d(v^{\a}\otimes v^{\g})=\sum_{i=1}^n \Omega_Q(\a,\g;i) v^{\a+[i]}\otimes v^{\g-[i]},
\end{equation}
with coefficients:
\begin{equation} \label{omq}
\begin{array}{l}
\Omega_Q(\a,\g;i)=
(-1)^{\ssum_{k<i}\g_k}{\disy \left(\prod_{k<i}q_{k,i}^{\g_k} \right) \left(\prod_{k>i}q_{k,i}^{\a_k} \right)\times}
\qquad \qquad \\
\hfill
{\disy \left( 1-p_i\left(\prod_{k=1}^n q_{i,k}^{\a_k+\g_k}\right)\right) }
 \ \textrm{ if } \g_i=1, \\
\textrm{and}\\
\Omega_Q(\a,\g;i)=0\ \textrm{ if } \g_i=0.
\end{array}
\end{equation}

\begin{sprop}[\cite{Si}, Proposition 3.5] \label{acy}
The complex
 $(K^{\sigma}(S_QV),d)$ above computes the twisted
 Ho\-ch\-schild homology of $S_{Q}V$, and ${\rm deg}(v^{\a}\otimes v^{\g})=\a+\g$ defines a   ${\mathbb N}^n$-grading on it.
Moreover, set
  $C^{\sigma}(Q)=\{ \rho\in {\mathbb N}^n \ |\ \forall i, \ \rho_i=0 \textrm{ or }\  p_iv_iv^{\rho}=v^{\rho}v_i \}$. Then for all
 $\rho \in {\mathbb N}^n\setminus C^{\sigma}(Q)$, the homogeneous subcomplex of $(K(S_QV),d)$ of degree
 $\rho$ is acyclic.
Further, the twisted Hochschild homology of $S_Q(V)$ in degree $k$ is given by
$$H_k(S_QV,{}_{\sigma}(S_QV))=\bigoplus_{\tiny{\begin{array}{c}\gamma\in\{0,1\}^n\\ |\gamma|=k\end{array}}}\bigoplus_{\tiny{\begin{array}{c}\alpha\in{\mathbb N}^n\\ \alpha+\gamma\in C^{\sigma}(Q)\end{array}}}\CC.v^{\alpha}\otimes v^{\gamma}.$$
\end{sprop}

The proof of this result relies on a homotopy $h_Q$ given by the
linear map $h_Q:K^{\sigma}_*\to K^{\sigma}_{*+1}$ defined by
$h_Q(v^{\a}\otimes
v^{\beta})={1\over{||\a+\b||}}\sum_{k=1}^n\omega_Q(\a,\b,i)v^{\a-\epsilon_i}\otimes
v^{\b+\epsilon_i}$, with
\begin{equation}\label{petitomega}
\omega_Q(\a,\b,i)=\left\{\begin{array}{ll} 0& {\rm if}\ \a+\b\in C^{\sigma}(Q)\\
0& {\rm if}\ \b_i=1\\
0& {\rm if}\ \a_i=0\\
\Omega_Q(\a-\epsilon_i,\b+\epsilon_i,i)^{-1}&{\rm otherwise}
\end{array}\right.
\end{equation}
where $||\gamma||$ is the cardinal of $\{k\ |\ \sigma(v_k)v_{\gamma}\neq v_{\gamma}v_k \mbox{ and } \gamma_k\neq 0 \}$, see the proof of \cite[Theorem 6.1]{Wa} and \cite[Proposition 3.5]{Si}.

\medskip

From the free resolution
 $U\otimes \L^*_QV\otimes U$ of $U$ we derive a complex $(R^*,{}^t\partial)$ which computes the Hochschild cohomology of $U$ with values in a bimodule $M$.
As a  $\CC$-vector space,
$$R^*=\Hom_{U^e}(U\otimes \L^*_QV\otimes U,M),$$
and the  differential is the  transposition of the  differential
$$\partial :U\otimes \L_Q^{*+1}V\otimes U\to U\otimes \L_Q^{*}V\otimes U,$$
 that is:
$$\begin{array}{rl}
{}^t\partial : & R^* \to R^{*+1}\\
 & \varphi \mapsto \varphi\circ \partial.
\end{array}$$
Then
$HH^*(U)=H^*(R,{}^t\partial)$.

\subsection{Duality} \label{gendu}

\hspace{1cm}
We present in this section  the links  between the Hochschild homology and  cohomology of the algebra $S_QV$. This section is mainly inspired from Section 6 of \cite{Rhh}, although the notation we use is slightly different.
We recall the following easy result, which will be our main tool in the following.

\begin{slemme}[\cite{Rhh}, Lemma 6.2.1] \label{trsp}
Let $(C_*, d)$ be a $\CC$-differential complex, and let $M_*$ be a graded $\CC$-vector space, such that there exists an isomorphism $\Phi$ of graded vector spaces  of degree 0, with source $C_*$ and target $M_*$.
Then the map $\wt d=\Phi\circ d\circ \Phi^{-1}$ is such that ${\wt d}^2=0$, and $(M_*, \wt d)$ is a differential complex. The  map $\Phi$  is then an isomorphism of complexes, and one has:
$$H_*(C,d)=H_*(M,\wt d).$$
\end{slemme}


Let us apply this result to the complexes
 $(R^*,{}^t\partial)$ and  $(K^{\sigma}_*,d)$ described  above.
We just give formulae here, the proofs can be found in Section 6 of \cite{Rhh}. Note that the quantum affine space we consider here is just a particular case of the ``mixed crossed algebras'' studied in \cite{Rhh}, and that all the following could apply {\sl verbatim} to these algebras.

There is an isomorphism $\Phi_{1,*}$ from $\Hom_{U^e}(U\otimes \L^* V\otimes U,U)$ onto $\Hom_k(\L^*_QV,U)$ defined for all $\varphi\in\Hom_{U^e}(U\otimes \L^* V\otimes U,U)$ by:
$$\Phi_{1,*}(\varphi)(v_{i_1}\wedge\ldots\wedge v_{i_{*}})= \varphi(1\otimes v_{i_1}\wedge\ldots\wedge v_{i_{*}}\otimes 1).$$
Then one computes the conjugate
 $D=\Phi_{1,*+1}\circ{}^t\partial\circ\Phi_{1,*}^{-1}$ \hbox{ of } ${}^t\partial$ \hbox{ by } $\Phi_1$.
Set
$$L^*=\Hom_k(\L^*_QV,U).$$
Then one has
$$\begin{array}{rrl}
D :& L^* & \to L^{*+1}\\
& \varphi &\mapsto D(\varphi), \end{array}$$
where $D(\varphi)$ is defined by:
\begin{equation} \label{d1}
\begin{array}{l}
D(\varphi)(v_{i_1}\wedge\ldots\wedge v_{i_{*+1}})=\\
\noalign{\medskip}
\blanc \blanc{\disy \sum_{k=1}^{*+1}(-1)^{k-1}\Big( (\prod_{s<k}q_{i_s,i_k})v_{i_k}\varphi(v_{i_1}\wedge\ldots \what v_{i_k}\ldots\wedge v_{i_{*+1}})}\\
\noalign{\medskip}
 \hfill - {\disy (\prod_{s>k}q_{i_k,i_s})\varphi(v_{i_1}\wedge\ldots \what v_{i_k}\ldots\wedge v_{i_{*+1}})v_{i_k}\Big).}\end{array}
\end{equation}

By construction, the complex $(R^*,{}^t\partial)$  and the complex $(L^*,D)$ above are isomorphic, and the  diagram below commutes:

$$\xymatrixrowsep{4.4pc}\xymatrixcolsep{6pc}\def\labelstyle{\large}
\xymatrix{R^* \ar[r]^{{}^t\partial} \ar[d]^{\Phi_{1,*}} \ar@{}[dr]|{\circlearrowleft} &R^{*+1} \ar[d]^{\Phi_{1,*+1}}\\
\Hom(\L_Q^*V,U)\ar[r]^{D}&\Hom(\L_Q^{*+1}V,U)
}
$$

\smallskip

Next, we prove a similar result for the complex $(K^{\sigma}_*,d)$, the homology of which is the Hochschild homology of $U$ with values in ${}_{\sigma}U$.
Recall that $K^{\sigma}_*= U \otimes\L^*_QV$ as a vector space and that the differential $d$ is given by formula (\ref{dwq}).

\begin{slemme}[\cite{Rhh}, Lemma 6.2.3]
The canonical map $\psi_* :\L^*_QV\otimes\L_Q^{n-*}V\to k\otimes v_1\wedge\ldots\wedge v_{n}$ defined by $\psi_*(v_{i_1}\wedge\ldots\wedge v_{i_*}\otimes v_{j_1}\wedge\ldots\wedge v_{j_{n-*}})=v_{i_1}\wedge\ldots\wedge v_{i_*}\wedge v_{j_1}\wedge\ldots\wedge v_{j_{n-*}}$ induces an isomorphism $\surl\psi_*: \L_Q^{n-*}V\to (\L^*_QV)'$ defined by $$\surl\psi_*(v_{j_1}\wedge\ldots\wedge v_{j_{n-*}})=\psi_*(\cdot\otimes v_{j_1}\wedge\ldots\wedge v_{j_{n-*}}).$$
\end{slemme}

\begin{sremark}{\rm
In fact $\surl\psi_*$ is just the linear map sending the element $v_{j_1}\wedge\ldots\wedge v_{j_{n-*}}$  to $\Theta_*(i_1,\ldots,i_*)(v_{i_1}\wedge\ldots\wedge v_{i_*})'$ where $\{i_1,\ldots,i_*\}$ is the complementary $*$-uple of $\{j_1,\ldots,j_{n-*}\}$ and $\Theta_*(i_1,\ldots,i_*)\in {\mathbb C}^*$ is defined by:
\begin{equation} \label{ll}
\Theta_*(i_1,\ldots,i_*) =
{\disy \prod_{k<i_*,\ k\not\in\{i_s\}} (-q_{i_*,k}) \prod_{k<i_{*-1},\ k\not\in\{i_s\}} (-q_{i_{*-1},k})\ldots \prod_{k<i_1,\ k\not\in\{i_s\}} (-q_{i_1,k}).}
\end{equation}}
\end{sremark}

The isomorphism $\surl\psi_*$ induces an isomorphism $\Phi_{2, *}={\rm id}\otimes\surl\psi_*$ from ${}_{\sigma}U\otimes\L_Q^{n-*}V$ to the space ${}_{\sigma}U\otimes(\L_Q^*V)' $.
But ${}_{\sigma}U\otimes\L_Q^{n-*}V=K^{\sigma}_{n-*}$, and one thus defines a differential $\Delta$ on the complex ${}_{\sigma}U\otimes(\L_QV)' $, such that the following diagram commutes:

$$\xymatrixrowsep{4.4pc}\xymatrixcolsep{6.8pc}\def\labelstyle{\large}
\xymatrix{{}_{\sigma}U\otimes(\L_Q^*V)'\ar@{}[rd]|\circlearrowleft \ar[r]^{\Delta} &{}_{\sigma}U\otimes(\L_Q^{*+1})'\\
K^{\sigma}_{n-*}\ar[u]^{\Phi_{2,*}}\ar[r]^d &K^{\sigma}_{n-*-1}\ar[u]^{\Phi_{2,*+1}}}
$$

The differential $\Delta$ is exactly $\Phi_{2,*+1}\circ d\circ\Phi_{2,*}^{-1}$.
For all  $i_1<\ldots<i_*$,
let $\{j_1,\ldots,j_{n-*}\}=\surl{\{i_1,\ldots,i_*\}}$  be the complementary set.
Then
\begin{equation} \label{d2}
\begin{array}{l}
\Delta(a\otimes(v_{i_1}\wedge\ldots\wedge v_{i_*})')=\Theta_*^{-1}(i_1,\ldots,i_*)\times \\
\blanc{\disy \sum_{k=1}^{n-*}(-1)^{k-1}\Big((\prod_{s<k}q_{j_s,j_k})av_{j_k}-(\prod_{s>k}q_{j_k,j_s})p_{j_k}v_{j_k}a\Big)}\\
\blanc \blanc \otimes \Theta_{*+1}(i_1,\ldots,j_k,\ldots,i_*)(v_{i_1}\wedge\ldots v_{j_k}\ldots \wedge v_{i_*})'.
\end{array}
\end{equation}

Once again, by construction the homology of the above complex
is:
$$H^*(U\otimes(\L_QV)')=HH_{n-*}(U,{}_{\sigma}U).$$

We have transfered the differential complex structures of $(K^{\sigma}_*,d)$ and  $(R^*,{}^t\partial)$ onto the graded vector spaces $\Hom(\L^*_QV,U)$ and  $U\otimes (\L^*_QV)'$.
There is a natural linear isomorphism between these two spaces;
we use it to compare the two differential complex structures.

Let $\Phi_{3,*}$ be the isomorphism  from $U\otimes (\L^*_QV)'$ to $\Hom(\L^*_QV,U)$ defined by:
\begin{equation} \label{f3}
\Phi_{3,*}(a\otimes\varphi)(v_{i_1}\wedge\ldots\wedge v_{i_*})=\varphi(v_{i_1}\wedge\ldots\wedge v_{i_*})a.
\end{equation}
Then consider the diagram below:

$$\xymatrixrowsep{4.4pc}\xymatrixcolsep{6pc}\def\labelstyle{\large}
\xymatrix{\Hom(\L_Q^*V,U)\ar[r]^D&\Hom(\L_Q^{*+1}V,U)\\
{}_{\sigma}U\otimes(\L_Q^*V)'\ar[u]^{\Phi_{3,*}}\ar[r]^{\Delta}&{}_{\sigma}U\otimes(\L_Q^{*+1}V)'\ar[u]^{\Phi_{3,*+1}}
}
$$
This diagram does not commute {\sl a priori}, but we have the following result.
The following Proposition relies directly on the computations leading to Lemma 6.3.2 of \cite{Rhh}. But we first must note that a careful check of these computations show that the formulation of this Lemma is wrong in the following sense. In expressing $\omega'_2$, the $\tilde\lambda_{j_k,t}$ must be replaced by $\tilde\lambda_{t,j_k}$. This does not affect the rest of \cite{Rhh}, since the only case considered there is the one where the products among $t$ of these elements is always equal to 1. Once this is observed, one gets the following.
\begin{sprop}
With the notation above, we have:
 $\Phi_{3,*+1}\circ \Delta=(-1)^{*+1}D \circ \Phi_{3,*}$.
\end{sprop}
\begin{demo}
 Lemmas 6.3.1 and 6.3.2 from \cite{Rhh}  can be rewritten in our context in the following form (taking care that the Hochschild homology is twisted):
$\Phi_{3,*+1}\circ \Delta (a\otimes(v_{i_1}\wedge\ldots\wedge v_{i_*})')$ is the linear map which sends the basis element
 $v_{\a_1}\wedge\ldots\wedge v_{\a_{*+1}}\in\L_Q^{*+1}V$ to:
$$\sum_{k=1}^{n-*}(-1)^{k-1} (\omega_1(\a_1,\ldots,\a_{*+1};k)av_{j_k}-p_{j_k}\omega_2(\a_1,\ldots,\a_{*+1};k)v_{j_k}a),
$$
and
$D\circ \Phi_{3,*} (a\otimes(v_{i_1}\wedge\ldots\wedge v_{i_*})')$ is the linear map which sends the basis element $v_{\a_1}\wedge\ldots\wedge v_{\a_{*+1}}\in\L_Q^{*+1}V$ to:
\begin{equation} \label{pli}
\sum_{k=1}^{n+r-*}(-1)^{k-1} (\omega'_1(\a_1,\ldots,\a_{*+1};k)av_{j_k}-\omega'_2(\a_1,\ldots,\a_{*+1};k)v_{j_k}a),
\end{equation}
with
$\omega'_1=(-1)^{*+1}\omega_1$, and
 $\omega'_2=(-1)^{*+1}p_{j_k}\omega_2$; so we are done.
\end{demo}

\medskip

The results of this section can now be gathered together in the following

\begin{sprop}
Let $U=S_QV$ be a quantum affine space. Then there is a duality between the Hochschild cohomology of $U$ and its Hochschild homology with values in the bimodule ${}_{\sigma}U$:
$$H_*(U,{}_{\sigma}U)\equiv H^{n-*}(U,U).$$
\end{sprop}

This result can be seen as a corollary of Van den  Bergh's Theorem
\cite{VdB}. We present the proof above in order to show how the
duality occurs at the level of the complexes themselves.

\begin{sremark}{\rm
Let $\L$ be an $m\times m$ multiplicatively skew-symmetric matrix. Define $Q=\left(\begin{array}{cc} \L&{}^t\L\\{}^t\L&\L\end{array}\right)$ by block, an  $n\times n$ multiplicatively skew-symmetric matrix with $n=2m$.
Then obviously the canonical automorphism $\sigma$ associated to $Q$ is the identity.
The duality result above implies in particular that for such a matrix $Q$ the Hochschild homology in degree n is nonzero, so that there is no dimension drop of the Hochschild homology with value in the algebra itself. This may explain a result by Connes and Dubois-Violette in the setting of smooth functions over a quantum real affine space ${\mathbb R}^{2m}$ parametrised by such a matrix, see \cite[Theorem 8]{connesdbviolette}.}
\end{sremark}

\section{Twisted Poincar\'e duality for the semiclassical limit of a quantum affine space}

\hspace{1cm} Now, we consider the uniparametrised case, by which we mean that the entries $q_{ij}$ of the matrix $Q$ are all powers of a generic $q$. Denote $q_{ij}=q^{a_{ij}}$, with $a_{ij}\in {\mathbb Z}$. As $q$ is generic, we have $a_{ii}=a_{ij}+a_{ji}=0$ for all $i,j$.
The semiclassical limit, \cite[Section III.5.4]{bg}, of the quantum
affine space $S_QV$ is the commutative algebra $R=\CC[X_1,\ldots,X_n]$
endowed with the Poisson bracket defined by
$\{X_i,X_j\}=a_{ij}X_iX_j$.

 The canonical Poisson homology and cohomology are defined
respectively thanks  to the K\"ahler differentials and the
multiderivations of $R$. We will denote the homology and
cohomology spaces by $HP_*(R)$ and $HP^*(R)$. The complexes computing
these homology groups are explicitly written down in the case where $n=2$ for instance
in \cite{RV} (see also \cite{Pichexp}), where the Poisson cohomology
is computed for the affine plane for any Poisson structure defined by
a homogenous polynomial.

In this section  we proceed as follows. First, we
consider the semiclassical limit $M$ of the twisted bimodule structure of
the  quantum space ${}_{\sigma}({S}_QV)$. It turns out that $M$ is
a Poisson dualising $R$-module, in the sense that there is a twisted
Poincar\'e duality between the Poisson cohomology of $R$ and its
Poisson homology with values in $M$:
\begin{eqnarray}
\label{dualitypoisson}
HP_k(R,M) \simeq HP^{n-k}(R).
\end{eqnarray}
 Next, we compute the Poisson homology of $R$ with values in $M$. As a
 consequence, we compute the Poisson cohomology of $R$ thanks to the
 above isomorphism (\ref{dualitypoisson}), and so retrieve a result of Monnier, \cite{Monnier}.

\subsection{Poisson algebra and Poisson module.}
\label{defM}

\hspace{1cm}
A commutative algebra $R$ endowed with a Lie bracket $\{.,.\}$ such that, for all $r \in R$, the map
$\{r,.\} : R \rightarrow R$ is a $\CC$-linear derivation of $R$ is called a Poisson algebra. From \cite{maszczyk,OH}, a Poisson module over the Poisson algebra $R$ is a $\CC$-vector space $M$ endowed with two bilinear maps $.$ and $\{.,.\}_M$ such that
\begin{enumerate}
\item $(M,.)$ is a (right) module over the commutative algebra $R$,
\item $(M,\{.,.\}_M)$ is a (right) module over the Lie algebra $(R, \{.,.\})$,
\item $x.\{a,b\} = \{x,a\}_M.b-\{x.b,a\}_M$ for all $a,b \in R$ and $x \in M$.
\item $\{x,ab\}_M=\{x,a\}_M.b +\{x,b\}_M.a$ for all $a,b \in R$ and $x \in M$.
\end{enumerate}

Starting from $S_QV$, we will now exhibit the Poisson structure on the
polynomial algebra $R=\CC[X_1,\dots,X_n]$ that arises from the semiclassical
limit process, see \cite[Section III.5.4]{bg}. Let  $A=(a_{ij}) \in
\mathcal{M}_n(\mathbb{Z})$ be antisymmetric and set $q_{ij}:=q^{a_{ij}}$, where $q$ is generic.
The semiclassical limit, \cite[Section III.5.4]{bg} of the quantum
affine space $S_QV$ is the commutative algebra $R=\CC[X_1,\ldots,X_n]$
endowed with the Poisson bracket defined on the generators of $R$ as
follows:
$$\{X_i,X_j\}:=([v_i,v_j]/(q-1))_{|q=1}=((1-q^{a_{ji}})/(q-1))_{|q=1}X_iX_j=a_{ij}X_iX_j.$$

Our aim in this section is to construct a Poisson module over the
Poisson algebra $R$ (with the Poisson structure as above) that will
restore a duality between Poisson homology and Poisson cohomology. As the dualising bimodule is ${}_{\sigma}({S}_QV)$ in the quantum setting, we will now consider the semiclassical limit
of this bimodule. In Section \ref{sectionduality}, we will show that the Poisson module resulting
from this procedure is actually a ``dualising Poisson module''.

We denote by $M$ the semiclassical limit of the twisted bimodule structure of the  quantum space ${}_{\sigma}({S}_QV)$.
As a vector space, $M={\mathbb C}[X_1,\ldots,X_n]=R$, and $M$ is
endowed with the following two actions of $R$:
\begin{itemize}
\item
the external product ``$.$'' is just the usual product of $R$.
\item the external bracket $\{.,.\}_M$ is defined by
$$\{m, X_i\}_M:=\left.\frac{mv_i-\sigma(v_i)m}{q-1}\right|_{q=1}$$
for all $m \in M$ and $i \in \{1,\dots , n\}$. In particular, when $m
= X_1^{\a_1}\ldots X_n^{\alpha_n}$ is a monomial
$$\{X_1^{\a_1}\ldots X_n^{\alpha_n}, X_i\}_M =\ \left.\frac{v^{\alpha}v_i-\sigma(v_i)v^{\a}}{q-1}\right|_{q=1}\ =\ \left.\frac{(\prod_{j>i}q_{ji}^{\a_j}-\prod_jq_{ji}\prod_{j<i}q_{ij}^{\a_j})v^{\a+\epsilon_i}}{q-1}\right|_{q=1}.$$
Recall that $q_{ij}=q^{a_{ij}}$;
so $$(\prod_{j>i}q_{ji}^{\a_j}-\prod_jq_{ji}\prod_{j<i}q_{ij}^{\a_j})=\prod_{j>i}q_{ji}^{\a_j}(1-q^{\sum_j a_{ij}(\a_j-1)}),$$
 and finally $\{X_1^{\a_1}\ldots X_n^{\alpha_n}, X_i\}_M=-\sum_j a_{ij}(\a_j-1) X^{\a+\epsilon_i}$.
\end{itemize}

One can easily check that $M$ is a Poisson module over $R$. Further, observe
 that

\begin{eqnarray}
\label{twistmodule}
 \{m, X_i\}_M = -\left\{X_i,m \right\} + \left(\sum_{l=1}^n a_{il} \right)
 X_im,
\end{eqnarray}
 for all $m \in M$.

\subsection{Poisson homology}

\hspace{1cm}
Let $M$ be a Poisson module over a Poisson algebra $R$. Then one defines a chain complex on the $R$-module
$C_{*}^{Poiss}(R,M)=\oplus_{k\in \mathbb{N}} C_k^{Poiss}(R,M)$, where $C_k^{Poiss}(R,M):= M\otimes_R \Omega^k(R) $ and $\Omega^k(R)$ denotes the so-called K\"ahler differential  $k$-forms, as follows, \cite{maszczyk}.
The boundary operator $\partial_k:  C_k^{Poiss}(R,M) \rightarrow  C_{k-1}^{Poiss}(R,M)$ is defined by
$$\begin{array}{l}
\hspace{-.8cm} \partial_k (m \otimes da_1 \wedge \dots \wedge da_k )  =  \ssum_{i=1}^k (-1)^{i+1} \{m ,a_i\}_M \otimes da_1 \wedge \dots \wedge \widehat{da_i} \wedge \dots \wedge da_k \\
 \blanc \blanc \esp + \ssum_{1 \leq i < j \leq k} (-1)^{i+j} m \otimes d\{a_i,a_j\} \wedge da_1 \wedge \dots \wedge \widehat{da_i} \wedge \dots \wedge \widehat{da_j} \wedge \dots \wedge da_k ,
\end{array}$$
where we have removed the expressions under the hats in the previous sums and $d$ denotes the exterior differential.

One can easily check that $\partial_k $ is well-defined and that
$\partial_{k-1} \circ \partial_{k}=0$. The homology of this
complex is denoted by $HP_{*}(R,M)$. This homology is called the
canonical homology in \cite{maszczyk} in reference to the
canonical homology defined by Brylinski, \cite{brylinski}. It is also
called the Poisson homology of the Poisson algebra $R$ with
values in the Poisson module $M$.

In the particular case that we will study, $R$ will be a (commutative) polynomial algebra over $\CC$, namely $R=\CC[X_1, \dots, X_n]$. In this case, it is clear that $\Omega^*(R)$ is the $R$-module generated by the wedge products of the $1$-differential forms $dX_1, \dots , dX_n$.

\subsection{Poisson cohomology}

\hspace{1cm}
We denote by $\chi^k(R)$ the $R$-module of all skew-symmetric $k$-linear derivations of $R$, that is, the set of all skew-symmetric $\CC$-linear maps $R^k \rightarrow R$ that are derivations in each of their arguments. Then we set $\chi^*(R):= \oplus_{k \in \mathbb{N}} \chi^k(R)$, the $R$-module of so-called skew-symmetric multiderivations of $R$. One can define a cochain complex structure on this $R$-module as follows. The Poisson coboundary operator $\delta_k : \chi^k(R) \rightarrow \chi^{k+1}(R)$ is defined by
\begin{eqnarray*}
\delta_k(P)(f_0,\dots,f_k) & := &  \sum_{i=0}^k(-1)^i \left\{ f_i, P(f_0,\dots,\widehat{f_i},\dots,f_k)\right\}\\
& &  +\sum_{0 \leq i < j \leq k} (-1)^{i+j}P\left( \{f_i,f_j\},f_0, \dots,\widehat{f_i},\dots,\widehat{f_j}, \dots, f_k \right)
\end{eqnarray*}
for all $P \in \chi^k(R)$. It is easy to check that $\delta_k(P)$ belongs indeed to $\chi^{k+1}(R)$ and that
$\delta_{k+1} \circ \delta_k =0$. The cohomology of this complex is called the Poisson cohomology of $R$; it is denoted by $HP^*(R)$.

\subsection{A Poincar\'e duality result.}
\label{sectionduality}

\hspace{1cm}
For the rest of this paper, we assume that  $R=\CC[X_1,\dots,X_n]$ is  endowed with the Poisson bracket defined by
$$\{X_i,X_j\}=a_{ij}X_iX_j,$$
where $A=(a_{ij}) \in \mathcal{M}_n(\mathbb{Z})$ is skew-symmetric.

Recall  that, thanks to the canonical volume form $dX_1\wedge\ldots\wedge dX_n$, the set $\chi^k(R)$ of all
skew-symmetric $k$-linear derivations of $R$ is isomorphic as a vector
space to $\Omega^{n-k}(R) $ via an isomorphism $\dag$ defined as
follows. We denote by $S_n$ the set of all $n$-permutations. Further,
for all $\sigma \in S_n$, we denote by $\varepsilon(\sigma)$ its sign
and we set $\sigma_i:=\sigma(i)$. For all $P \in \chi^k(R)$ let $\dag(P)$ be the unique element of  $\Omega^{n-k}(R) $ defined by
$$\dag(P)=\sum_{\sigma \in S_{k,n-k}}\epsilon(\sigma) P(X_{\sigma_1}, \dots, X_{\sigma_k}) dX_{\sigma_{k+1}} \wedge \dots \wedge dX_{\sigma_n},$$
where $S_{k,n-k}$ denotes the set of those permutations $\sigma \in S_n$ such that $\sigma_1 < \dots < \sigma_k$ and
$\sigma_{k+1} < \dots < \sigma_{n}$.

From now on, $M$ denotes the Poisson $R$-module defined in Section \ref{defM}. Recall that  $M=\CC[X_1,\dots,X_n]=R$ as a vector space.
Hence, we deduce from the above result that  $\chi^k(R)$ is isomorphic as a vector space to $M \otimes_R \Omega^{n-k}(R) $ via an isomorphism still denoted by $\dag$ and defined by:

$$\dag(P)=\sum_{\sigma \in S_{k,n-k}}\epsilon(\sigma) P(X_{\sigma_1}, \dots, X_{\sigma_k}) dX_{\sigma_{k+1}} \wedge \dots \wedge dX_{\sigma_n}$$
for all $P \in \chi^k(R)$. (Observe that we have omitted $\otimes$ in order to simplify  the notation.)

So we have a diagram as follows.
$$\xymatrixrowsep{4.4pc}\xymatrixcolsep{6pc}\def\labelstyle{\large}
\xymatrix{\chi^k(R)\ar[d]^{\delta_k}\ar[r]^{\dag}&M \otimes_R \Omega^{n-k}(R)\ar[d]^{\partial_{n-k}}\\
\chi^{k+1}(R)\ar[r]^{\dag}&M \otimes_R \Omega^{n-k-1}(R)
}
$$

In order to prove that there is a (twisted) Poincar\'e duality between
the Poisson homology of $R$ with values in $M$ and the cohomology of
$R$, we will prove that this diagram is almost commutative.

\begin{sprop} \label{propositiondualite}
For all $P \in \chi^k(R)$, the following equality holds:
$$(\dag \circ \delta)(P)=(-1)^{k+1} (\partial_{} \circ \dag)(P).$$
\end{sprop}
\begin{demo} First, it follows from the definition of $\delta$ and $\dag$ that  $(\dag \circ \delta)(P)  =  U +V$, where
$$U:=\zum{\sigma \in S_{k+1,n-k-1}}{1\leq i\leq k+1}\varepsilon(\sigma) (-1)^{i+1} \left\{ X_{\sigma_i}, P[X_{\sigma_1}, \dots , \widehat{X_{\sigma_i}}, \dots, X_{\sigma_{k+1}}] \right\} dX_{\sigma_{k+2}} \wedge \dots \wedge dX_{\sigma_n} $$
and
$$ V=  \zum{\sigma \in S_{k+1,n-k-1}}{1\leq i<j\leq k+1}\varepsilon(\sigma) (-1)^{i+j}  P[ \{X_{\sigma_i},X_{\sigma_j}\}, X_{\sigma_1}, \dots , \widehat{X_{\sigma_i}}, \dots,\widehat{X_{\sigma_j}}, \dots, X_{\sigma_{k+1}}]
dX_{\sigma_{k+2}} \wedge \dots \wedge dX_{\sigma_n}.$$

We proceed in three steps.\\

\noindent $\bullet$ {\bf Step 1: we rewrite $U$.}\\

For $r$ distinct integers $i_1,\ldots,i_r$ let $(i_1\ldots i_r)$ be the cyclic permutation sending $i_1$ to $i_2$, $\ldots$, $i_{r-1}$ to $i_r$ and $i_r$ to $i_1$.
Let $S_{k,1,n-k-1}$ be the set of those permutation $\tau \in S_n$ such that $\tau_1 < \dots < \tau_k$ and $\tau_{k+2} < \dots < \tau_n$. Then,  the map $(\sigma,i) \in S_{k+1,n-k-1} \times  \{1, \dots ,k+1\} \mapsto \tau \in S_{k,1,n-k-1}$ given by
$$\tau = \sigma \circ (i \: i+1 \dots k+1) $$ is well-defined and
turns out to be a bijection. Indeed, it is easy to see that these two
sets have each cardinality $\left(\!{\tiny\begin{array}{c}n\\ k\end{array}}\!\right)\times(n-k)$ and that this map is injective. Observe
further that $\varepsilon(\tau)=(-1)^{k+1-i}
\varepsilon(\sigma)$. Thus, by means of the change of variables
induced by this bijection, we get:
\begin{eqnarray}
\label{expressionU}
U =  (-1)^{k}\sum_{\tau \in S_{k,1,n-k-1}} \varepsilon(\tau) \left\{
X_{\tau_{k+1}}, P[X_{\tau_1}, \dots , X_{\tau_k}] \right\}
dX_{\tau_{k+2}} \wedge \dots \wedge dX_{\tau_n}.\end{eqnarray}
$ $

\noindent $\bullet$ {\bf Step 2: we rewrite $V$.}\\

As $ \{X_{\sigma_i},X_{\sigma_j}\}=a_{\sigma_i,\sigma_j}  X_{\sigma_i}X_{\sigma_j}$ and $P$ is a skew-symmetric multiderivation, one can rewrite $V$ as follows.

\begin{eqnarray*}
\hspace{-15mm}V  =&  & \hspace{-5mm} \zum{\sigma \in S_{k+1,n-k-1}}{1 \leq i < j \leq k+1}\varepsilon(\sigma)  (-1)^{i+j} a_{\sigma_i,\sigma_j} (-1)^jX_{\sigma_i} P[ X_{\sigma_1}, \dots , \widehat{X_{\sigma_i}}, \dots, X_{\sigma_{k+1}}]  dX_{\sigma_{k+2}} \wedge \dots \wedge dX_{\sigma_n} \\
 & + &\hspace{-5mm}\zum{\sigma \in S_{k+1,n-k-1}}{1 \leq i < j \leq k+1}\varepsilon(\sigma)  (-1)^{i+j}a_{\sigma_i,\sigma_j} (-1)^{i+1} X_{\sigma_j} P[  X_{\sigma_1},  \dots,\widehat{X_{\sigma_j}}, \dots, X_{\sigma_{k+1}}]  dX_{\sigma_{k+2}} \wedge \dots \wedge dX_{\sigma_n}\\
 = & & \hspace{-5mm}\sum_{\sigma \in S_{k+1,n-k-1}}\varepsilon(\sigma) \sum_{1 \leq i < j \leq k+1} (-1)^{i} a_{\sigma_i,\sigma_j} X_{\sigma_i} P[ X_{\sigma_1}, \dots , \widehat{X_{\sigma_i}}, \dots, X_{\sigma_{k+1}}]  dX_{\sigma_{k+2}} \wedge \dots \wedge dX_{\sigma_n} \\
 & + &\hspace{-5mm} \sum_{\sigma \in S_{k+1,n-k-1}}\varepsilon(\sigma) \sum_{1 \leq i < j \leq k+1} (-1)^{j+1}a_{\sigma_i,\sigma_j} X_{\sigma_j} P[  X_{\sigma_1},  \dots,\widehat{X_{\sigma_j}}, \dots, X_{\sigma_{k+1}}]  dX_{\sigma_{k+2}} \wedge \dots \wedge dX_{\sigma_n}\\
 = & &\hspace{-5mm} \sum_{\sigma \in S_{k+1,n-k-1}}\varepsilon(\sigma) \sum_{i=1}^k (-1)^{i}\left(\sum_{j=i+1}^{k+1} a_{\sigma_i,\sigma_j}\right) X_{\sigma_i} P[ X_{\sigma_1}, \dots , \widehat{X_{\sigma_i}}, \dots, X_{\sigma_{k+1}}]  dX_{\sigma_{k+2}} \wedge \dots \wedge dX_{\sigma_n} \\
& +&\hspace{-5mm} \sum_{\sigma \in S_{k+1,n-k-1}}\varepsilon(\sigma) \sum_{j=2}^{k+1} (-1)^{j} \left(\sum_{i=1}^{j-1}a_{\sigma_j,\sigma_i}\right) X_{\sigma_j} P[  X_{\sigma_1},  \dots,\widehat{X_{\sigma_j}}, \dots, X_{\sigma_{k+1}}]  dX_{\sigma_{k+2}} \wedge \dots \wedge dX_{\sigma_n}\\
=& &\hspace{-5mm} \sum_{\sigma \in S_{k+1,n-k-1}}\varepsilon(\sigma)
 \sum_{i=1}^{k+1} (-1)^{i}\left(\sum_{j=1}^{k+1}
 a_{\sigma_i,\sigma_j}\right) X_{\sigma_i} P[ X_{\sigma_1}, \dots ,
 \widehat{X_{\sigma_i}}, \dots, X_{\sigma_{k+1}}]  dX_{\sigma_{k+2}}
 \wedge \dots \wedge dX_{\sigma_n}
\end{eqnarray*}
Hence, using the same change of variables as in the previous step, we
get:
\begin{eqnarray}
\label{expressionV}
V =  (-1)^{k+1}\sum_{\tau \in S_{k,1,n-k-1}}\varepsilon(\tau)  \left(\sum_{l=1}^{k} a_{\tau_{k+1},\tau_l}\right) X_{\tau_{k+1}} P[ X_{\tau_1}, \dots , X_{\tau_{k}}]  dX_{\tau_{k+2}} \wedge \dots \wedge dX_{\tau_n}
\end{eqnarray}
$ $

\noindent $\bullet$ {\bf Step 3: we rewrite  $(\partial \circ \dag)(P)$ and conclude.} \\

First, it follows from the definition of $\partial$ and $\dag$ that
$$\begin{array}{l}(\partial \circ \dag)(P)   =   \ssum_{\sigma \in S_{k,n-k}} \varepsilon(\sigma) \left[ \ssum_{i=1}^{n-k}
(-1)^{i+1} \big\{P[X_{\sigma_1}, \dots,X_{\sigma_k}],X_{\sigma_{i+k}}\big\}_M  dX_{\sigma_{k+1}} \wedge \dots \wedge \widehat{dX_{\sigma_{k+i}}} \wedge \dots \wedge dX_{\sigma_n} \right. \\
 \left. + \ssum_{1 \leq i < j \leq n-k}
(-1)^{i+j} P[X_{\sigma_1}, \dots,X_{\sigma_k}]  d\left\{X_{\sigma_{k+i}}, X_{\sigma_{k+j}} \right\} \wedge dX_{\sigma_{k+1}} \wedge \dots \wedge \widehat{dX_{\sigma_{k+i}}} \wedge \dots \wedge \widehat{dX_{\sigma_{k+j}}}\dots \wedge dX_{\sigma_n}\right].
\end{array}$$
Next, using (\ref{twistmodule}) and the fact that $\left\{X_{\sigma_{k+i}}, X_{\sigma_{k+j}} \right\}=a_{\sigma_{k+i}\sigma_{k+j}} X_{\sigma_{k+i}} X_{\sigma_{k+j}}$, we obtain:
\begin{eqnarray*}
 (\partial \circ \dag)(P) &= & \sum_{\sigma \in S_{k,n-k}}\varepsilon(\sigma) \sum_{i=1}^{n-k}
(-1)^i \left\{X_{\sigma_{i+k}},P[X_{\sigma_1}, \dots,X_{\sigma_k}]\right\} dX_{\sigma_{k+1}} \wedge \dots \wedge \widehat{dX_{\sigma_{k+i}}} \wedge \dots \wedge dX_{\sigma_n} \\
& & \hspace{-3cm}+   \sum_{\sigma \in S_{k,n-k}}\varepsilon(\sigma) \sum_{i=1}^{n-k}
(-1)^{i+1} \left(\sum_{l=1}^n a_{\sigma_{i+k} l}\right) X_{\sigma_{i+k}} P[X_{\sigma_1}, \dots,X_{\sigma_k}] dX_{\sigma_{k+1}} \wedge \dots \wedge \widehat{dX_{\sigma_{k+i}}} \wedge \dots \wedge dX_{\sigma_n} \\
& & \hspace{-3cm} + \sum_{\sigma \in S_{k,n-k}}\varepsilon(\sigma) \sum_{1 \leq i < j \leq n-k}
(-1)^{i+j} a_{\sigma_{k+i} \sigma_{k+j}}(-1)^{j-2}X_{\sigma_{k+i}} P[X_{\sigma_1}, \dots,X_{\sigma_k}]  dX_{\sigma_{k+1}} \wedge \dots \wedge \widehat{dX_{\sigma_{k+i}}} \wedge \dots \wedge dX_{\sigma_n} \\
& & \hspace{-3cm} + \sum_{\sigma \in S_{k,n-k}}\varepsilon(\sigma) \sum_{1 \leq i < j \leq n-k}
(-1)^{i+j} a_{\sigma_{k+i} \sigma_{k+j}}(-1)^{i-1}X_{\sigma_{k+j}} P[X_{\sigma_1}, \dots,X_{\sigma_k}]  dX_{\sigma_{k+1}} \wedge \dots \wedge \widehat{dX_{\sigma_{k+j}}} \wedge \dots \wedge dX_{\sigma_n}
\end{eqnarray*}
Then, rewriting the last three sums in the right-hand side leads to
\begin{eqnarray*}
(\partial \circ \dag)(P) &  = & \sum_{\sigma \in S_{k,n-k}}\varepsilon(\sigma) \sum_{i=1}^{n-k}
(-1)^i \left\{X_{\sigma_{i+k}},P[X_{\sigma_1}, \dots,X_{\sigma_k}]\right\} dX_{\sigma_{k+1}} \wedge \dots \wedge \widehat{dX_{\sigma_{k+i}}} \wedge \dots \wedge dX_{\sigma_n} \\
& & \hspace{-2cm} +   \sum_{\sigma \in S_{k,n-k}}\varepsilon(\sigma) \sum_{i=1}^{n-k}
(-1)^{i+1} \left(\sum_{l=1}^k a_{\sigma_{i+k} \sigma_l}\right) X_{\sigma_{i+k}} P[X_{\sigma_1}, \dots,X_{\sigma_k}] dX_{\sigma_{k+1}} \wedge \dots \wedge \widehat{dX_{\sigma_{k+i}}} \wedge \dots \wedge dX_{\sigma_n}
\end{eqnarray*}

Finally, observe that the map $(\sigma,i) \in S_{k,n-k} \times  \{1, \dots ,n-k\} \mapsto \tau \in S_{k,1,n-k-1}$ given by
$$\tau = \sigma \circ (k+i \dots k+1) $$ is well-defined and is a bijection. Observe further that
$\varepsilon(\tau)=(-1)^{i-1} \varepsilon(\sigma)$. Hence, by means of the change of variable induced by this bijection,  we get
\begin{eqnarray*}
 (\partial \circ \dag)(P)  &=  & -\sum_{\tau \in S_{k,1,n-k-1}}\varepsilon(\tau)
 \left\{X_{\tau_{k+1}},P[X_{\tau_1}, \dots,X_{\tau_k}]\right\} dX_{\tau_{k+2}} \wedge  \dots \wedge dX_{\tau_n} \\
& & \hspace{-1cm} +   \sum_{\tau \in S_{k,1,n-k-1}}\varepsilon(\tau)  \left(\sum_{l=1}^k a_{\tau_{k+1} \tau_l}\right) X_{\tau_{k+1}} P[X_{\tau_1}, \dots,X_{\tau_k}] dX_{\tau_{k+2}} \wedge  \dots \wedge dX_{\tau_n}
\end{eqnarray*}
So, we deduce from (\ref{expressionU}) and (\ref{expressionV}) that $$(\partial \circ \dag)(P) =(-1)^{k+1} (U+V)=(-1)^{k+1} (\dag \circ \delta_{\pi})(P), $$
as desired.
\end{demo}

${}$\\

It follows from Proposition \ref{propositiondualite} that the diagram
$$\xymatrixrowsep{4.4pc}\xymatrixcolsep{6pc}\def\labelstyle{\large}
\xymatrix{\chi^k(R)\ar[r]^{\dag}\ar[d]^{\delta_k}\ar@{}[rd]|\circlearrowleft
  & M \otimes_R \Omega^{n-k}(R)\ar[d]^{(-1)^{k+1}\partial_{n-k}}\\
\chi^{k+1}(R)\ar[r]^{\dag} & M\otimes_R \Omega^{n-k-1}(R)
}
$$
is commutative. Naturally, this leads to a (twisted) Poincar\'e duality between the Poisson homology of $R$ with values in $M$ and the Poisson cohomology of $R$. More precisely, we can now state the main result of this paper.

\begin{stheoreme}
\label{duality}
For all $k \in \mathbb{N}$, we have:
$$ HP_k(R,M) \simeq HP^{n-k}(R).$$
\end{stheoreme}

\subsection{Application to Poisson cohomology.}

\hspace{1cm}
In view of Theorem \ref{duality}, in order to compute the Poisson
cohomology of $R$, it is enough to compute $ HP_k(R,M)$. The final
part of this paper is dedicated to this computation.

First, if $\alpha=(\alpha_1,\dots,\alpha_n) \in \mathbb{N}^n$, then we set $X^{\alpha}:=X_1^{\alpha_1} \dots X_n^{\alpha_n}$. Similary,
if $\beta=(\beta_1,\dots,\beta_n) \in \mathbb{N}^n$, then we set
$dX^{\beta}:=dX_1^{\beta_1} \wedge \dots \wedge dX_n^{\beta_n}$. Then,
by following the proof of \cite[Theorem 6.1]{Wa}, we obtain the following result.

\begin{sprop}
\label{homology}
$$HP_k(R,M)=\ozum{\mid \beta \mid =k }{ \alpha +\beta \in C} \mathbb{C} X^{\alpha} dX^{\beta}$$
where
$$C:=\{ \gamma\in{\mathbb N}^n \mid \gamma_i=0 \mbox{ or } \sum_{j=1}^na_{ij}(\gamma_j-1)=0\}.$$
\end{sprop}

\begin{demo}  Easy computations show that
$$\partial(X^{\alpha}dX^{\beta} )
=\sum_{i=1}^{n} \delta_{1,\beta_i}\Omega(\alpha,\beta,i) X^{\alpha+\varepsilon_i} dX^{\beta-\varepsilon_i},$$
where $\delta_{1,\beta_i}$ is the Kronecker symbol, and
$$\Omega(\alpha,\beta,i):=(-1)^{\sum_{j=1}^{i-1}\beta_j} \sum_{j=1}^n a_{ij}(\alpha_j +\beta_j -1).$$

Observe that $\Omega(\alpha,\beta,i)
=\frac{\Omega_Q(\alpha,\beta,i)}{1-q}\mid_{q=1}$, where
$Q:=(q^{a_{ij}})$ and $\Omega_Q(\alpha,\beta,i)$ has been defined in (\ref{omq}).
We also set
$$\omega(\alpha,\beta,i):=\ (1-q)\omega_Q(\alpha,\beta,i) \mid_{q=1},$$
where $\omega_Q$ has been defined in (\ref{petitomega}).

It is clear that $X^{\alpha}dX^{\beta}$ is in the homology group of the complex when $\alpha + \beta \in C$.
Next we have to prove that there exists an homotopy $h$ such that $\partial_{\pi,\sigma}h+h\partial_{\pi,\sigma}(X^{\alpha}dX^{\beta})=X^{\alpha}dX^{\beta}$ when $\alpha + \beta \notin C$.
We set
$$h(X^{\alpha}dX^{\beta}):= \frac{1}{\mid \mid \alpha + \beta \mid  \mid }\sum_{i=1}^n\omega(\alpha,\beta,i)X^{\alpha-\varepsilon_i}dX^{\beta+\varepsilon_i}.$$

In order to prove that $h$ is a homotopy, it is enough to prove
that certain equalities hold between the $\Omega$ and the
$\omega$. As these equalities hold at  the quantum level, i.e. the
linear map $h_Q$ defined by (\ref{petitomega}) is a homotopy, the
desired equalities also hold at the ``semiclassical level'' thanks
to a specialisation.
\end{demo}
$ $\\

In the case where $n=2$ and $a_{12}=1$, we obtain the following result:
$$HP_0(R,M)={\mathbb C}\oplus {\mathbb C}X_1X_2,\ \ HP_1(R,M)={\mathbb C}X_1dX_2\oplus {\mathbb C}X_2dX_1,\ \ HP_2(R,M)={\mathbb C}dX_1\wedge dX_2.$$
In this way, we retrieve the dimensions computed for the cohomology in \cite{RV}, see also the Introduction of the present work.

Finally, in view of Theorem \ref{duality} and Proposition
\ref{homology}, we obtain the following result regarding the
Poisson cohomology of $R$. This result has been previously
obtained by a direct (and so completely different) method in
\cite{Monnier}.

\begin{scorol}

$$HP^k(R) \simeq \ozum{\mid \beta \mid =n-k }{ \alpha +\beta \in C} \mathbb{C} X^{\alpha} dX^{\beta}$$
where
$$C:=\{ \gamma\in{\mathbb N}^n \mid \gamma_i=0 \mbox{ or } \sum_{j=1}^na_{ij}(\gamma_j-1)=0\}.$$
\end{scorol}

\section*{Aknowledgements}
\hspace{1cm}
The authors thank Pol Vanhaecke for bringing the paper \cite{Monnier} to our attention, and Tom Hadfield for the reference \cite{Si}. We also thank Ulrich Kr\"ahmer for comments on a previous version of this paper.

Part of this work was done while the first author was visiting the Newton Institute. He wishes to thank Shahn Majid for his kind invitation to participate to the Noncommutative Geometry program.

The main ideas of this paper where initiated when both authors attended an algebra conference in La Rochelle in February 2006. We are grateful to the organisers and the participants of this conference.

It is our pleasure to thank Tom Lenagan for fruitful discussions on the topic of this paper, and for comments on an earlier version.

\bigskip

\bigskip

\noindent
St\'ephane Launois, Lionel Richard\\
University of Edinburgh, School of Mathematics and Maxwell Institute for Mathematical Sciences \\
 JCMB - Mayfield Road\\
  Edinburgh EH9 3JZ, United Kingdom.

  \smallskip

  \noindent
   stephane.launois@ed.ac.uk,  \ \ lionel.richard@ed.ac.uk

\end{document}